\documentclass[a4paper,10pt]{article}
\usepackage[utf8]{inputenc}

\title{RKFD Methods - a short review}
\author{Maciej Jaromin\footnote{University of Wroc\l{}aw}}
\date{\today}

\usepackage{indentfirst}
\usepackage[utf8]{inputenc}
\usepackage[sumlimits,intlimits]{amsmath}
\usepackage{mathrsfs,amssymb,amsthm,mathtools,bbm, amsthm}

\usepackage[T1]{fontenc}
\usepackage{lmodern,hyperref}

\usepackage{geometry}
\theoremstyle{definition}
\makeatletter
\newtheorem{przykl@d}[theorem]{Example}
\newenvironment{example}{\begin{przykl@d}}{\qed\end{przykl@d}}
\makeatother
\theoremstyle{remark}

\begin{document}
\maketitle

\begin{abstract}
In this paper, a recently published method [Hussain, Ismail, Senua, Solving directly special fourth-order ordinary differential equations using Runge–Kutta type method, J. Comput. Appl. Math. 306 (2016) 179–199] for solving fourth-order ordinary differential equations is summarized and reviewed. An independant implementation of one of the published methods is also tested and benchmarked against the RK4 method.
\end{abstract}

\section{Introduction} \label{intro}

In a recent paper \cite{Name}, the authors proposed a direct numerical method of solving special cases of the initial value problem (IVP) for fourth-order ordinary differential equations (ODEs) called the RKFD method. Solving fourth-order ODEs has applications in multiple fields, including beam theory \cite{uses1, uses2}, fluid dynamics \cite{uses3}, neural networks \cite{uses4}, electric circuits \cite{uses5} and the ill-posed problem of a beam on elastic foundation \cite{uses6}.

The fourth-order initial value problem is a differential equation of the form
\[y^{(\textrm{\romannumeral 4})}\left(x\right) = f\left(x,y(x),y'(x),y''(x),y'''(x)\right),\, x \geq x_0\] with the initial conditions $y(x_0) = y_0$, $y'(x_0) = y'_0$, $y''(x_0) = y''_0$, $y'''(x_0) = y'''_0$. However, the RKFD method is only applicable in the special case where $y^{(\textrm{\romannumeral 4})}\left(x\right) = f\left(x,y(x)\right)$ -- that is, it does not contain the first, second or third derivative of $y$.

Although fourth-order ODEs can be solved by transforming them into an equivalent system of first-order equations (see \cite[p. 283]{system1}, \cite[p. 4]{system2}, \cite[p. 253]{system3}), such methods are computationally expensive. Other direct methods of solving fourth-order ODEs have been proposed before \cite{multistep1,multistep2,multistep3,multistep4,multistep5,multistep6}. However, the RKFD method is a one-step method, while the ones cited above are two-step methods. This gives RKFD an advantage, as it requires less function evaluations to compute, leading to better efficiency.

This paper is organized as follows: Section 2 gives an overview of how the general RKFD method is formulated in the original paper. Section 3 covers the algebraic order conditions for the RKFD method and their use to construct explicit RKFD methods. Section 4 describes the numerical results obtained by the authors of the original paper and presents the results of implementing the explicit fourth-order RKFD4 method and comparing it to the standard fourth-order Runge-Kutta method (RK4).

\section{RKFD method formulation}
In this section, the process of formulating the RKFD method is summarized. The general RKFD method is formulated \cite[\S 2]{Name} by converting the fourth order IVP into a system of first order ODEs:
\[\left(\begin{array}{c}y(x)\\ v(x)\\ u(x)\\ w(x)\end{array}\right)' = \left(\begin{array}{c}v(x)\\ u(x)\\ w(x)\\ f(x, y(x))\end{array}\right)\]
The general Runge-Kutta method is then applied, with the parameters $A = \{a_{ij}\}$, $b \{b_{i}\}$, $c = \{c_{i}\}$ obtaining a system of equations \cite[eqs. (4)-(11)]{Name} which is the starting point for the formulation of the RKFD method. Because $f$ only contains $x$ and $y(x)$, the system can be simplified significantly \cite[eqs. (12)-(16)]{Name}. The equations are then further simplified using some well-known properties of consistent Runge-Kutta methods, such as: $$\sum_{j=1}^{s} a_{ij} = c_i,\; \text{for } i = 1, ..., s$$
These operations define the RKFD method, which takes the following form \cite[eqs. (17)-(21)]{Name}:
\begin{align}
y_{n+1} &= y_n + hy'_n + \frac{1}{2}h^2y''_n +\frac{1}{6}h^3y'''_n + h^4\sum_{i=1}^s b_if (x_n + c_ih, Y_i),\\
y_{n+1}' &= y'_n + hy''_n + \frac{1}{2}h^2y'''_n +h^3\sum_{i=1}^s b'_if (x_n + c_ih, Y_i),\\
y_{n+1}'' &= y''_n + hy'''_n +h^2\sum_{i=1}^s b''_if (x_n + c_ih, Y_i),\\
y_{n+1}''' &= y'''_n + h\sum_{i=1}^s b'''_if (x_n + c_ih, Y_i),\\
Y_i &= y_n +  c_ihy'_n + \frac{1}{2}c_i^2h^2y''_n +\frac{1}{6}c_i^3h^3y'''_n + h^4 \sum_{j=1}^s \hat{a}_{ij}f (x_n + c_jh, Y_j),
\label{eq:1}
\end{align}
where $b' = \{b'_{i}\}$, $b'' = \{b_{i}''\}$, $b''' = \{b_{i}'''\}$ and $\hat{A} = \{\hat{a}_{ij}\}$ are defined as such in relation to the parameters $A$, $b$, $c$ of the Runge-Kutta method:
$$ b^TA = b'^T,\, b^TA^2 = b''^T,\, b^TA^3 = b'''^T,\, A^4 = \hat{A}$$
The parameters of a RKFD method can be written in a manner similar to the Butcher tableau of a RK method:
\begin{table}[h]
\centering
\begin{tabular}{c|c}
		$c$ & $\hat{A}$ \\
		\hline\\[-2ex]
		& $b^T$ \\
		& $b'^T$ \\
		& $b''^T$ \\
		& $b'''^T$ \\
\end{tabular}
=
\begin{tabular}{c|ccc}
	$c_1$ & $\hat{a}_{11}$ & $\hdots$ & $\hat{a}_{1s}$ \\
	$\vdots$ & $\vdots$ & $\ddots$ & $\vdots$ \\
	$c_s$ & $\hat{a}_{s1}$ & $\hdots$ & $\hat{a}_{ss}$ \\ 
	\hline\\[-2ex]
	& $b_1$ & $\hdots$ & $b_s$\\
	& $b'_1$ & $\hdots$ & $b'_s$\\
	& $b''_1$ & $\hdots$ & $b''_s$\\
	& $b'''_1$ & $\hdots$ & $b'''_s$\\
\end{tabular}
\caption{The parameters of a general RKFD method written in the form of a tableau.}
\label{table:1}
\end{table}

It is worth noting that in \eqref{eq:1}, the authors of the original paper have made a typographical error \cite[eq. (21)]{Name}. A similar error is also present earlier in the paper when presenting the general s-stage Runge-Kutta method \cite[eq. (3)]{Name}. In both cases, the sum in the definition of $Y_i$ is defined over both $i$ and $j$, instead of just $j$.

\section{Constructing instances of the RKFD method} 

The authors of \cite{Name} present a number of algebraic conditions required for an RKFD method to have a specific order. These conditions are obtained by representing both the exact solution of the problem and the numerical solution obtained through the RKFD method as B-series. The order conditions for the RKFD method are as follows (see \cite[eqs. (73)-(79)]{Name}):

Order 1:
\begin{equation}
b'''^Te = 1,
\end{equation}

Order 2:
\begin{equation}
b'''^Tc = \frac{1}{2},\, b''^Te = \frac{1}{2},
\end{equation}

Order 3:
\begin{equation}
b'''^Tc^2 = \frac{1}{3},\, b''^Tc = \frac{1}{6},\, b'^Te = \frac{1}{6},
\end{equation}

Order 4:
\begin{equation}
b'''^Tc^3 = \frac{1}{4},\, b''^Tc^2 = \frac{1}{12},\, b'^Tc = \frac{1}{24},\, b^Te = \frac{1}{24},
\end{equation}

Order 5:
\begin{equation}
b'''^Tc^4 = \frac{1}{5},\, b'''^T\hat{A} = \frac{1}{120},\, b''^Tc^3 = \frac{1}{20},\, b'^Tc^2 = \frac{1}{60},\, b^Tc = \frac{1}{120},
\end{equation}

Order 6:
\begin{equation}
b'''^Tc^5 = \frac{1}{6},\, b'''^T\hat{A}c = \frac{1}{720},\, b'''^T(c.\hat{A}e) = \frac{1}{144},\, b''^Tc^4 = \frac{1}{30},\, b''^T\hat{A} = \frac{1}{720},\, b'^Tc^3 = \frac{1}{120},\, b^Tc^2 = \frac{1}{360},
\end{equation}

Order 7:
\begin{equation}
\begin{aligned}
b'''^Tc^6 = \frac{1}{7},\, b'''^T(c.\hat{A}c) = \frac{1}{840},\, b'''^T(c^2.\hat{A}e) = \frac{1}{168},\, b'''^T(\hat{A}c^2) = \frac{1}{2520},\, b''^Tc^5 = \frac{1}{42}, \\b''^T\hat{A}c = \frac{1}{5040}, b''^T(c.\hat{A}e) = \frac{1}{1008},\, b'^Tc^4 = \frac{1}{210},\, b'^T\hat{A} = \frac{1}{5040},\, b^Tc^3 = \frac{1}{840},
\end{aligned}
\end{equation}

$e$ is not explicitly defined in \cite{Name}, but by analyzing \cite[tables 2 and 3]{Name}, it can be determined that $e$ is a vector of ones.

To construct an explicit RKFD method of a particular order, the order requirements are used to obtain a system of equations. While the system does not have a single solution, some of the parameters are selected as free variables and their values are subsequently obtained through the method of minimizing error equations (a technique estabilished by Dormand and Prince \cite{minimize}). Using this process, minimizing error equations using Maple software, two explicit three-stage RKFD methods are constructed, a fourth-order method designated RKFD4 and a fifth-order method designated RKFD5 (see tables \ref{table:2}, \ref{table:3}).
\begin{table}[h]
	\centering
\begin{tabular}{c|ccc}
	$\frac{4}{11}$ & $-\frac{1}{5}$ \\[6pt]
	$\frac{17}{20}$ & $\frac{19}{125}$  & $\frac{19}{125}$\\[6pt]
	\hline\\[-2ex]
	& $\frac{17}{200}$ & $-\frac{7}{75}$ & $\frac{1}{20}$\\[6pt]
	& $\frac{1}{18}$ & $\frac{209}{1926}$ & $\frac{6}{1926}$\\[6pt]
	& $\frac{47}{408}$ & $\frac{847}{2568}$ & $\frac{100}{1819}$\\[6pt]
	& $\frac{47}{408}$ & $\frac{1331}{2568}$ & $\frac{2000}{5457}$\\[6pt]
\end{tabular}
\caption{The Butcher tableau for the RKFD4 method.}
\label{table:2}

\begin{tabular}{c|ccc}
	$\frac{3}{5}+\frac{\sqrt{6}}{10}$ & $\frac{4059}{187793}$ \\[6pt]
	$\frac{3}{5}-\frac{\sqrt{6}}{10}$ & $-\frac{1502}{532215}$  & $\frac{1826}{569317}$\\[6pt]
	\hline\\[-2ex]
	& $\frac{19}{1080}$ & $\frac{13}{1080}-\frac{11\sqrt{6}}{2160}$ & $\frac{13}{1080}+\frac{11\sqrt{6}}{2160}$\\[6pt]
	& $\frac{1}{18}$ & $\frac{1}{18}-\frac{\sqrt{6}}{48}$ & $\frac{1}{18}+\frac{\sqrt{6}}{48}$\\[6pt]
	& $\frac{1}{9}$ & $\frac{7}{36}-\frac{\sqrt{6}}{18}$ & $\frac{7}{36}+\frac{\sqrt{6}}{18}$\\[6pt]
	& $\frac{1}{9}$ & $\frac{4}{9}-\frac{\sqrt{6}}{36}$ & $\frac{4}{9}+\frac{\sqrt{6}}{36}$\\[6pt]
\end{tabular}
\caption{The Butcher tableau for the RKFD5 method.}
\label{table:3}
\end{table}

\section{Numerical results and implementation}
In this section, the numerical results obtained in the original paper are described, along with the results of implementating the RKFD4 method and comparing its performance with the RK4 method. The RKFD4 and RKFD5 methods were tested on the following set of problems in \cite[\S 6]{Name}:\\
\begin{example}
	Problem 1:\\
	$y^{(\textrm{\romannumeral 4})} = -4y,$\\
	$y(0)=0,\,y'(0)=1,\,y''(0)=2,\,y'''(0)=2,\,$\\
	Integrated over the interval $[0,10]$. Exact solution: $y(x)=e^x\sin(x)$\\\\
\end{example}
\begin{example} \label{ex}
	Problem 2:\\
$y^{(\textrm{\romannumeral 4})} = y^2 + \cos^2(x) + \sin(x) - 1,$\\
$y(0)=0,\,y'(0)=1,\,y''(0)=0,\,y'''(0)=-1,\,$\\
Integrated over the interval $[0,10]$. Exact solution: $y(x)=\sin(x)$\\\\
\end{example}
\begin{example}
	Problem 3:\\
$y^{(\textrm{\romannumeral 4})} = \frac{3 \sin(y) (3+2\sin^2(y))}{cos^7(y)},$\\
$y(0)=0,\,y'(0)=1,\,y''(0)=0,\,y'''(0)=1,\,$\\
Integrated over the interval $[0,\frac{\pi}{4}]$. Exact solution: $y(x)=\arcsin(x)$\\\\
\end{example}
\begin{example}
	Problem 4:\\
$y^{(\textrm{\romannumeral 4})} = e^{3x}u, \,y(0)=1,\,y'(0)=-1,\,y''(0)=1,\,y'''(0)=-1,\,$\\
$z^{(\textrm{\romannumeral 4})} = 16e^{-x}y,\, z(0)=1,\,z'(0)=-2,\,z''(0)=4,\,z'''(0)=-8,\,$\\
$w^{(\textrm{\romannumeral 4})} = 81e^{-x}z,\, w(0)=1,\,w'(0)=-3,\,w''(0)=9,\,w'''(0)=-27,\,$\\
$u^{(\textrm{\romannumeral 4})} = 256e^{-x}w,\, u(0)=1,\,u'(0)=-4,\,u''(0)=16,\,u'''(0)=-64,\,$\\
Integrated over the interval $[0,2]$.\\
Exact solution: $y(x)=e^{-x},\,z(x)=e^{-2x},\,w(x)=e^{-3x},\,u(x)=e^{-4x}$\\\\
\end{example}
\begin{example}
	Problem 5: The ill-posed Problem of a Beam on Elastic Foundation.\\
$y^{(\textrm{\romannumeral 4})} = -y+1, 0<x<1,$\\
$y(0)=0,\,y'(0)=0,\,y''(0)=0,\,y'''(0)=0,\,$\\
Exact solution: $y(x)=1 - \frac{1}{2}e^{-\frac{x}{\sqrt{2}}}(1+e^{\sqrt{2x}})\cos(\frac{x}{\sqrt{2}})$
\end{example}

The performance of RKFD4 and RKFD5 was compared against RK4 (the classic four-stage fourth-order RK method given in Butcher \cite[p. 180]{butcher}), RK5Ns6 (a six-stage fifth-order RK method given in Butcher \cite[p. 192]{butcher}), RKN4D (a four-stage fourth-order Runge-Kutta-Nystram (RKN) method given in Dormand \cite[p. 265]{system3}) and RKN5H (a four-stage fifth-order RKN method given in Hairer \cite[p. 285]{system1}). In all five problems, the RKFD methods were more efficient than the RK and RKN methods they were compared against, with RKFD5 in particular achieving a smaller maximum error with less total function evaluations than any other method tested at each step size and for every problem (see \cite[Figs. 3-7]{Name}).

For the purposes of this review, the RKFD4 method was implemented in the Julia programming language (available at \url{https://github.com/MaciejJaromin/RKFD-method-implementation}). Example \ref{ex} was implemented and used to confirm that RKFD4 is a fourth-order method by calculating the error of the method in a single step at various step sizes.

The implementation was also used to benchmark the performance of RKFD4 against the performance of RK4 on the same task. However, the results were inconsistent with the ones presented in \cite{Name}. While RKFD4 achieved greater accuracy, its calculation time was approximately $150\%$ that of RK4 for all step sizes. The results are summarized in table \ref{table:4}:
\begin{table}[h]
\centering
\begin{tabular}{c|cr|cr}
	&RK4&&RKFD4\\
	\hline
	Step size&Error&Time&Error&Time\\
	\hline
	$0.1$&$7.66\mathrm{e}{-04}$&0.017 s&$6.09\mathrm{e}{-04}$&0.028 s\\		$0.01$&$7.78\mathrm{e}{-08}$&0.204 s&$7.38\mathrm{e}{-09}$&0.309 s\\
	$0.001$&$7.78\mathrm{e}{-12}$&1.967 s&$2.11\mathrm{e}{-13}$&3.168 s\\
	$0.0001$&$7.78\mathrm{e}{-16}$&19.579 s&$1.59\mathrm{e}{-17}$&31.283 s\\
	$0.00001$&$7.78\mathrm{e}{-20}$&195.778 s&$1.54\mathrm{e}{-21}$&315.901 s\\
\end{tabular}
\caption{Benchmarking results for the RK4 and RKFD4 methods using Example \ref{ex}.}
\label{table:4}
\end{table}

A possible cause for the slower performance of RKFD4 is the fact that Julia is a high-level language and the implementation of RK4 used for the benchmark uses matrix operations which may be more optimized. 

\section{Conclusions}

RKFD methods are a promising take on the problem of integrating fourth-order ODEs directly. However, this new method has several problems, the key and unaddressable one being that it can only solve specific kinds of fourth-order differential equations (see \S \ref{intro}). Regardless, higher-order RKFD methods or implicit RKFD methods offer potential new avenues of research and RKFD methods in general offer a specialized but powerful numerical integration tool.

\end{document}